\documentclass[11pt]{article}
\usepackage{amssymb}
\newtheorem{lemma}{Lemma}[section]
\newtheorem{definition}[lemma]{Definition}

\newtheorem{theorem}[lemma]{Theorem}
\newtheorem{proposition}[lemma]{Proposition}
\newtheorem{corollary}[lemma]{Corollary}
\newtheorem{example}[lemma]{Example}

\newtheorem{note}[lemma]{Remark}

\def\endproof{\hfill$\Box$}

\def\endproof{\hfill$\Box$}

\title{On structures in hypergraphs of models of a theory\footnote{This research was partially
supported by Committee of Science in Education and Science
Ministry of the Republic of Kazakhstan (Grant No. AP05132546) and
Russian Foundation for Basic Researches (Project No.
17-01-00531-a).}}

\author{}
\author{B.Sh.~Kulpeshov, S.V.~Sudoplatov}
\date{}

\begin{document}

\maketitle

\begin{abstract}
We define and study structural properties of hypergraphs of models
of a theory including lattice ones. Characterizations for the
lattice properties of hypergraphs of models of a theory, as well
as for structures on sets of isomorphism types of models of a
theory, are given.
\end{abstract}

{\bf Keywords:} hypergraph of models, elementary theory,
elementarily substructural set, lattice structure.

\bigskip
Hypergraphs of  models of a theory are derived objects allowing to
obtain an essential structural information about both given
theories and related semantic objects including graph ones
\cite{CCMCT14, Su013, Su08, Baik, SudKar16, KulSud17, Sud17,
KulSud171, KS12018}. Studying of hypergraphs of models of a theory
is closely related with a series of papers on description of
lattices of substructures \cite{Paris1, Gaifman, Paris2, Wilkie,
Schmerl1, Mills, Schmerl11, Schmerl2, Schmerl3, KosSch, Schmerl4,
Schmerl5, AQ}.

In the presented paper we define and study structural properties
of hypergraphs of models of a theory including lattice ones.
Characterizations for the lattice properties of hypergraphs of
models of a theory as well as for structures on sets of
isomorphism types of models of a theory are given.

\section{Preliminaries}

Recall that a {\it hypergraph\/} is a pair of sets $(X,Y)$, where
$Y$ is some subset of the Boolean $\mathcal{P}(X)$ of the set $X$.

Let $\mathcal{M}$ be some model of a complete theory $T$.
Following \cite{SudKar16}, we denote by $H(\mathcal{M})$ a family
of all subsets $N$ of the universe $M$ of $\mathcal{M}$ that are
universes of elementary submodels $\mathcal{N}$ of the model
$\mathcal{M}$: $H(\mathcal{M})=\{N\mid
\mathcal{N}\preccurlyeq\mathcal{M}\}$. The pair
$(M,H(\mathcal{M}))$ is called the {\it hypergraph of elementary
submodels\/} of the model $\mathcal{M}$ and denoted by
$\mathcal{H}(\mathcal{M})$\index{$\mathcal{H}(\mathcal{M})$}.

%\smallskip
\begin{definition} {\rm \cite{KulSud171}. Let $\mathcal{M}$ be a model of a theory
$T$ with a hypergraph $\mathcal{H}=(M,H(\mathcal{M}))$ of
elementary submodels, $A$ be an infinite definable set in
$\mathcal{M}$, of arity $n$: $A\subseteq M^n$. The set $A$ is
called {\em $\mathcal{H}$-free} if for any infinite set
$A'\subseteq A$, $A'=A\cap Z^n$ for some $Z\in H(\mathcal{M})$
containing parameters for $A$. Two $\mathcal{H}$-free sets $A$ and
$B$ of arities $m$ and $n$ respectively are called {\em
$\mathcal{H}$-independent} if for any infinite $A'\subseteq A$ and
$B'\subseteq B$ there is $Z\in H(\mathcal{M})$ containing
parameters for $A$ and $B$ and such that $A'=A\cap Z^m$ and
$B'=B\cap Z^n$.}
\end{definition}

%\smallskip
Note the following properties \cite{KulSud171}.

%\smallskip
1. Any two tuples of a $\mathcal{H}$-free set $A$, whose distinct
tuples do not have common coordinates, have same type.

Indeed, if there are tuples $\bar{a},\bar{b}\in A$ with ${\rm
tp}(\bar{a})\ne{\rm tp}(\bar{b})$ then for some formula
$\varphi(\bar{x})$ the sets of solutions of that formula and of
the formula $\neg\varphi(\bar{x})$ divide the set $A$ into two
nonempty parts $A_1$ and $A_2$, where at least one part, say
$A_1$, is infinite. Taking $A_1$ for $A'$ we have $A_1=A\cap Z^n$
for appropriate $Z\in H(\mathcal{M})$ and~$n$. Then by the
condition for tuples in $A$ we have $A_2\cap Z^n=\emptyset$ that
is impossible since $Z$ is the universe of an elementary submodel
of $\mathcal{M}$.

Thus the formula $\varphi(\bar{x})$, defining $A$, implies some
complete type in $S^n(\emptyset)$, and if $A$ is
$\emptyset$-definable then $\varphi(\bar{x})$ is a principal
formula.

In particular, if the set $A$ is $\mathcal{H}$-free and
$A\subseteq M$, then the formula, defining $A$, implies some
complete type in $S^1(\emptyset)$.

%\smallskip
2. If $A\subseteq M$ is a $\mathcal{H}$-free set, then $A$ does
not have nontrivial definable subsets, with parameters in $A$,
i.e., subsets distinct to subsets defined by equalities and
inequalities with elements in $A$.

%\smallskip
Indeed, if $B\subset A$ is a nontrivial definable subset then $B$
is defined by a tuple $\bar{a}$ of parameters in $A$, forming a
{\em finite} set $A_0\subset A$, and $B$ is distinct to subsets of
$A_0$ and to $A\setminus C$, where $C\subseteq A_0$. Then removing
from $A$ a set $B\setminus A_0$ or $(A\setminus B)\setminus A_0$,
we obtain some $Z\in H(\mathcal{M})$ violating the satisfiability
for $B$ or its complement. It contradicts the condition that $Z$
is the universe of an elementary submode of $\mathcal{M}$.

%\smallskip
3. If $A$ and $B$ are two $\mathcal{H}$-independent sets, where
$A\cup B$ does not have distinct tuples with common coordinates,
then $A\cap B=\emptyset$.

Indeed, if $A\cap B$ contains a tuple $\bar{a}$, then, choosing
infinite sets $A'\subseteq A$ and $B'\subseteq B$ with $\bar{a}\in
A'$ and $\bar{a}\notin B'$, we obtain $\bar{a}\in A'=A\cap Z^n$
for appropriate $Z\in H(\mathcal{M})$ and $n$, as so $\bar{a}\in
B\cap Z^n=B'$. This contradiction means that $A\cap B=\emptyset$.

\begin{definition} {\rm \cite{KulSud17}. The {\em complete union}
of hypergraphs $(X_i,Y_i)$, $i\in I$, is the hypergraph
$\left(\bigcup\limits_{i\in I}X_i,Y\right)$, where
$Y=\left\{\bigcup\limits_{i\in I}Z_i\mid Z_i\in Y_i\right\}$. If
the sets $X_i$ are disjoint, the complete union is called {\em
disjoint} too.  If the set $X_i$ form a $\subseteq$-chain, then
the complete union is called {\em chain}.}
\end{definition}

%\smallskip
By Property 3 we have the following theorem on decomposition of
restrictions of hypergraphs $\mathcal{H}$, representable by unions
of families of $\mathcal{H}$-independent sets.

\begin{theorem}\label{dcu} {\rm \cite{KulSud171}.}
A restriction of hypergraph $\mathcal{H}=(M,H(\mathcal{M}))$ to a
union of a family of $\mathcal{H}$-free $\mathcal{H}$-independent
sets $A_i\subseteq M$ is represented as a disjoint complete union
of restrictions $\mathcal{H}_i$ of the hypergraph $\mathcal{H}$ to
the sets $A_i$.
\end{theorem}

Proof. Consider a family of $\mathcal{H}$-independent sets
$A_i\subseteq M$. By Property 3 these sets are disjoint, and using
the definition of $\mathcal{H}$-independence we immediately obtain
that the union of restrictions $\mathcal{H}_i$ of $\mathcal{H}$ to
the sets $A_i$ is complete.
\endproof

\vskip 3mm

Recall that a subset $A$ of a linearly ordered structure $M$ is
called {\it convex} if for any $a, b\in A$ and $c\in M$ whenever
$a<c<b$ we have $c\in A$. A {\it weakly o-minimal structure} is a
linearly ordered structure $M=\langle M,=,<,\ldots \rangle$ such
that any definable (with parameters) subset of the structure $M$
is a union of finitely many convex sets in $M$.

In the following definitions $M$ is a weakly o-minimal structure,
$A, B\subseteq M$, $M$ be $|A|^+$-saturated, $p,q\in S_1(A)$ be
non-algebraic types.
\begin{definition}\rm \cite{bbs1}
We say that $p$ is not {\it weakly orthogonal} to $q$
($p\not\perp^w  q$) if there exist an $A$-definable formula
$H(x,y)$, $\alpha \in p(M)$ and $\beta_1, \beta_2 \in q(M)$ such
that $\beta_1 \in H(M,\alpha)$ and $\beta_2 \not\in H(M,\alpha)$.
\end{definition}

\begin{definition} \rm \cite{k2003} We say that $p$ is not {\it quite orthogonal} to
 $q$ ($p\not\perp^q q$) if there exists an $A$-definable bijection $f: p(M)\to q(M)$.
 We say that a weakly o-minimal theory is {\it quite o-minimal} if the notions of
 weak and quite orthogonality of 1-types coincide.
\end{definition}

In the work \cite{KS} the countable spectrum for quite o-minimal
theories with non-maximal number of countable models has been
described:

\begin{theorem}\label{KS_apal} Let $T$ be a quite o-minimal theory with non-maximal number
of countable models. Then $T$ has exactly $3^k\cdot 6^s$ countable
models, where $k$ and $s$ are natural numbers. Moreover, for any
$k,s\in\omega$ there exists a quite o-minimal theory $T$ having
exactly $3^k\cdot 6^s$ countable models.
\end{theorem}

Realizations of these theories with a finite number of countable
models are natural ge\-ne\-ra\-li\-za\-tions of Ehrenfeucht
examples obtained by expansions of dense linear orderings by a
countable set of constants, and they are called theories of  {\em
Ehrenfeucht type}. Moreover, these realizations are representative
examples for hypergraphs of prime models \cite{CCMCT14, Su08,
SudKar16}. We consider operators for hypergraphs allowing on one
hand to describe the decomposition of hypergraphs of prime models
for quite o-minimal theories with few countable models, and on the
other hand pointing out constructions leading to the building of
required hypergraphs by some simplest ones.

Having nontrivial structures like structures with some orders it
is assumed that ``complete'' decompositions are considered modulo
additional conditions guaranteing the elementarity for
substructures with considered universes. So we use the {\em
conditional} completeness taking unions with the properties of
density, linearity etc.

Below we illustrate this conditional completeness for structures
with dense linear orders.

Denote by $(M,H_{\rm dlo}(\mathcal{M}))$ the hypergraph of (prime)
elementary submodels of a countable model $\mathcal{M}$ of the
theory of dense linear order without endpoints.

\begin{note}\label{no441} {\rm  The class of hypergraphs
$(M,H_{\rm dlo}(\mathcal{M}))$ is closed under countable chain
complete unions, modulo density and having an encompassing dense
linear order without endpoints. Thus, any hypergraph $(M,H_{\rm
dlo}(\mathcal{M}))$ is represented as a countable chain complete,
modulo density, union of some its proper subhypergraphs.}
\end{note}

Any countable model of a theory of Ehrenfeucht type is a disjoint
union of some intervals, which are ordered both themselves and
between them, and of some singletons. Dense subsets of the
intervals form universes of elementary substructures. So, in view
of Remark \ref{no441}, we have:

\begin{theorem}\label{KS_apr} {\rm \cite{KulSud17}.}
A hypergraph of prime models of a countable model of a theory of
Ehrenfeucht type is represented as a disjoint complete, modulo
density, union of some hypergraphs in the form $(M,H_{\rm
dlo}(\mathcal{M}))$ as well as singleton hypergraphs of the form
$(\{c\},\{\{c\}\})$.
\end{theorem}

\begin{note}\label{no442} {\rm Taking into consideration links
between sets of realizations of $1$-types, which are not weakly
orthogonal, as well as definable equivalence relations, the
construction for the proof of Theorem \ref{KS_apr} admits a
natural generalization for an arbitrary quite o-minimal theory
with few countable models. Here conditional complete unions should
be additionally {\em coordinated}, i.e., considering definable
bijections between sets of realizations of $1$-types, which are
not quite orthogonal.}
\end{note}

\section{Elementarily substructural sets}

Let $\mathcal{M}$ be a model of theory $T$, $(M,H(\mathcal{M}))$
be a hypergraph of elementary submodels of $\mathcal{M}$. The sets
$N\in H(\mathcal{M})$ are called {\em elementarily submodel} or
{\em elementarily substructural} in $\mathcal{M}$.

Elementarily substructural sets in $\mathcal{M}$ are characterized
by the following well-known  Tarski--Vaught Theorem, which is
called the Tarski--Vaught test.

\begin{theorem}\label{thTV}
Let $\mathcal{A}$ and $\mathcal{B}$ be structures in a language
$\Sigma$, $\mathcal{A}\subseteq\mathcal{B}$. The following are
equivalent:

{\rm (1)} $\mathcal{A}\preccurlyeq\mathcal{B}${\rm ;}

{\rm (2)} for any formula $\varphi(x_0,x_1,\ldots,x_n)$ in the
language $\Sigma$ and for any elements $a_1,\ldots,a_n\in A$, if
$\mathcal{B}\models\exists x_0\,\varphi(x_0,a_1,\ldots,a_n)$ then
there is an element $a_0\in A$ such that
$\mathcal{B}\models\varphi(a_0,a_1,\ldots,a_n)$.
\end{theorem}

\begin{corollary}\label{coess1}
A set $N\subseteq M$ is elementarily substructural in
$\mathcal{M}$ if and only if for any formula
$\varphi(x_0,x_1,\ldots,x_n)$ in the language
$\Sigma(\mathcal{M})$ and for any elements $a_1,\ldots,a_n\in N$,
if $\mathcal{M}\models\exists x_0\,\varphi(x_0,a_1$, $\ldots$,
$a_n)$ then there is an element $a_0\in N$ such that
$\mathcal{M}\models\varphi(a_0,a_1,\ldots,a_n)$.
\end{corollary}

\begin{proposition}\label{press1}
Let $A$ be a definable set in an $\omega_1$-saturated model
$\mathcal{M}$ of a countable complete theory $T$. Then exactly one
of the following conditions is satisfied:

$(1)$ $A$ is finite and contained in any elementarily
substructural set in $\mathcal{M}$;

$(2)$ $A$ is infinite and has infinitely many distinct
intersections with elementarily substructural sets in
$\mathcal{M}$, and all these intersections are infinite; these
intersections can be chosen forming an infinite chain/antichain by
inclusion.
\end{proposition}

Proof. If $|A|<\omega$ then $A$ is contained in ${\rm
acl}(\emptyset)$, and so it is contained in any elementary
submodel of $\mathcal{M}$.

If $A=\varphi(\mathcal{M},\bar{a})$ is infinite, we construct a
countable submodel $\mathcal{N}_0\prec\mathcal{M}$ containing
parameters in $\bar{a}$. Since $A$ is infinite, the set $A\cap
N_0$ is countable. By compactness, since $\mathcal{M}$ is
$\omega_1$-saturated, the set $A\setminus N_0$ is infinite. Adding
to $N_0$ new elements of $A$ we construct a countable model
$\mathcal{N}_1$ such that
$\mathcal{N}_0\prec\mathcal{N}_1\prec\mathcal{M}$. Continuing the
process we build an elementary chain of models $\mathcal{N}_k$,
$k\in\omega$, such that $\mathcal{N}_k\prec\mathcal{M}$ and $A\cap
N_k\subset A\cap N_{k+1}$, $k\in\omega$.

Constructing the required antichain of intersections $A\cap N$
with elementarily substructural sets $N$, it suffices to use
\cite[Theorem 2.10]{KS12018} allowing to separate disjoint finite
sets, whose elements do not belong to ${\rm acl}(\emptyset)$.
\endproof

\medskip
The arguments for the proof of Proposition \ref{press1} stay valid
for a countable saturated model $\mathcal{M}$. Thus, we have the
following

\begin{proposition}\label{press2}
Let $A$ be a definable set in a countable saturated model
$\mathcal{M}$ of a small theory $T$. Then exactly one of the
following conditions is satisfied:

$(1)$ $A$ is finite and contained in any elementarily
substructural set in $\mathcal{M}$;

$(2)$ $A$ is infinite and has infinitely many distinct
intersections with elementarily substructural sets in
$\mathcal{M}$, and all these intersections are infinite; these
intersections can be chosen forming an infinite chain/antichain by
inclusion.
\end{proposition}

The following example illustrates that if $\mathcal{M}$ is not
saturated then the conclusions of assertions \ref{press1} and
\ref{press2} can fail.

\begin{example}\label{exess1}
{\rm Let a set $A$ is defined by a unary predicate $P$ and
includes infinitely many language constants $c_i$, $i\in I$. Then
there is, in the language $\{P\}\cup\{c_i\mid i\in I\}$, a
structure $\mathcal{M}$ having only finite set $A_0$ of elements
in $A$, which are not interpreted by constants. Since elementarily
substructural sets $N$ take all constants, there are only finitely
many possibilities for intersections $A\cap N$.}
\end{example}

In view of aforesaid arguments it is interesting to describe
possible cardinalities both for sets $H(\mathcal{M})$ and their
restrictions $H(\mathcal{M})\upharpoonright
A\rightleftharpoons\{A\cap N\mid N\in H(\mathcal{M})\}$ on
definable sets $A\subseteq M$.

Since in Example \ref{exess1} intersections $A\cap N$, taking all
constants $c_i$, can include an arbitrary subset of $A_0$, then
for this example we have $|H(\mathcal{M})\upharpoonright
A|=2^{|A_0|}$. The same formula holds for infinite sets $A_0$, but
in such a case the set $H(\mathcal{M})\upharpoonright A$ is
transformed from finite one directly to a set with continuum many
elements.

Note that for $\mathcal{H}$-free sets $A\subseteq M$, {\em modulo}
${\rm acl}(\emptyset)$ (i.e., for sets $A$, whose each subset
$B\subseteq A\setminus {\rm acl}(\emptyset)$ has a representation
$B\cup ({\rm acl}(\emptyset)\cap A)=A\cap N$ for some $N\in
H(\mathcal{M})$), the equality $|H(\mathcal{M})\upharpoonright
A|=2^{|A\setminus{\rm acl}(\emptyset)|}$ holds. Thus, we have the
following {\em dichotomy theorem}.

\begin{theorem}\label{thess2}
For any $\mathcal{H}$-free, modulo ${\rm acl}(\emptyset)$, set
$A\subseteq M$ its restriction to any elementary submodel
$\mathcal{M}_0\prec\mathcal{M}$ satisfies either
$|H(\mathcal{M}_0)\upharpoonright A|=2^n$ for some $n\in\omega$,
or $|H(\mathcal{M}_0)\upharpoonright A|=2^\lambda$ form some
$\lambda\geq\omega$.
\end{theorem}

Similar to Example \ref{exess1}, the following example illustrates
the dichotomy for hypergraphs of elementary submodels.

\begin{example}\label{exess2}
{\rm Consider the structure $\mathcal{M}$ of rational numbers,
$\langle{\bf Q},<,c_q\rangle_{q\in{\bf Q}}$, in which every
element is interpreted by a constant. This structure does not have
proper elementary substructures, therefore
$|H(\mathcal{M})|=1=2^0$. Extending $\mathcal{M}$ to a structure
$\mathcal{M}_1$ by addition of $n$ elements for pairwise distinct
$1$-types, defined by cuts, we have $|H(\mathcal{M}_1)|=2^n$. If
$\mathcal{M}$ is extended till a structure $\mathcal{M}_2$ by
addition of at least two elements of fixed cut or of infinitely
many elements for distinct cuts, then by density the summarized
number of added elements occurs infinite and
$|H(\mathcal{M}_2)|=2^\lambda$ holds for some
$\lambda\geq\omega$.}
\end{example}

At the same time there are examples of hypergraphs of elementary
submodels, for which the conclusion of Theorem \ref{thess2} fails.
For instance, as shown in \cite{Wilkie}, there are hypergraphs for
the theory of arithmetic of natural numbers such that
$|H(\mathcal{M})|=5$ and the lattice of elementary submodels is
isomorphic to the lattice $P_5$.

\section{Lattice structures associated with hypergraphs of models of a theory}

\medskip
For given structure $\mathcal{M}$ we define the structure
$L(\mathcal{M})=\langle H(\mathcal{M});\wedge,\vee\rangle$ by the
following relations for
$\mathcal{M}_1,\mathcal{M}_2\prec\mathcal{M}$:
$\mathcal{M}_1\wedge\mathcal{M}_2=\mathcal{M}_1\cap\mathcal{M}_2$
and $\mathcal{M}_1\vee\mathcal{M}_2=\mathcal{M}(M_1\cup M_2)$.

Consider the following question: when the structure
$L(\mathcal{M})$ is a lattice?

Clearly, answering this question we have to characterize the
conditions $\mathcal{M}_1\cap\mathcal{M}_2\prec\mathcal{M}$ and
$\mathcal{M}(M_1\cup M_2)\prec\mathcal{M}$. Assuming that
$\mathcal{M}$ is infinite, the structures
$\mathcal{M}_1\cap\mathcal{M}_2$ should be infinite too, in
particular, $M_1\cap M_2\ne\emptyset$. By \cite[Theorem
3.2]{SudKar16}, assuming that $\mathcal{M}$ is
$\lambda$-saturated, it can not contain separated sets $A$ and $B$
of cardinalities $<\lambda$, such that ${\rm acl}(A)\cap{\rm
acl}(B)=\emptyset$.

By Theorem \ref{thTV} we have the following theorems
characterizing the elementarity of sub\-struc\-t\-ures.

\begin{theorem}\label{thTV2}
Let $\mathcal{M}_1$ and $\mathcal{M}_2$ be elementary
substructures of structure $\mathcal{M}$ in a language $\Sigma$,
$M_1\cap M_2\ne\emptyset$. The following are equivalent:

{\rm (1)} $(\mathcal{M}_1\cap\mathcal{M}_2)\prec\mathcal{M}${\rm
;}

{\rm (2)} for any formula $\varphi(x_0,x_1,\ldots,x_n)$ of the
language $\Sigma$ and for any elements $a_1,\ldots,a_n\in M_1\cap
M_2$ if $\mathcal{M}\models\exists x_0\,\varphi(x_0,a_1, \ldots,
a_n)$ then there is an element $a_0\in M_1\cap M_2$ such that
$\mathcal{M}_i\models\varphi(a_0,a_1,\ldots,a_n)$, $i=1,2$.
\end{theorem}

\begin{theorem}\label{thTV3}
Let $\mathcal{M}_1$ and $\mathcal{M}_2$ be elementary
substructures of structure $\mathcal{M}$ in a language $\Sigma$.
The following are equivalent:

{\rm (1)} $\mathcal{M}(M_1\cup M_2)\prec\mathcal{M}${\rm ;}

{\rm (2)} for any formula $\varphi(x_0,x_1,\ldots,x_n)$ of the
language $\Sigma$ and for any elements $a_1,\ldots,a_n\in M_1\cap
M_2$ if $\mathcal{M}\models\exists x_0\,\varphi(x_0,a_1, \ldots,
a_n)$ then there is an element $a_0\in M(M_1\cup M_2)$ such that
$\mathcal{M}(M_1\cup M_2)\models\varphi(a_0,a_1,\ldots,a_n)$.
\end{theorem}

The following examples illustrate valuations of the conditions (2)
in Theorems \ref{thTV2} and \ref{thTV3}.

\begin{example}\label{exess3}
{\rm Consider a structure $\mathcal{M}$ in a graph language
$\{R^{(2)}\}$ with a symmetric irreflexive relation $R$ and
elements $a_1,a_2,a_3,a_4$ such that
$$R=\{[a_1,a_3],[a_1,a_4],[a_2,a_3],[a_2,a_4]\}.$$ The substructures
$\mathcal{M}_1\subset\mathcal{M}$ and
$\mathcal{M}_2\subset\mathcal{M}$ with the universes
$\{a_1,a_2,a_3\}$ and $\{a_1,a_2,a_4\}$ respectively satisfy the
formula $\varphi(a_1,a_2)\rightleftharpoons\exists
x(R(a_1,x)\wedge R(a_2,x))$ whereas
$\mathcal{M}_1\cap\mathcal{M}_2$ does not satisfy that formula
since appropriate elements for $x$ belong to $M_1\oplus M_2$.}
\end{example}

\begin{example}\label{exess4}
{\rm Consider a structure $\mathcal{M}$ of graph language
$\{R^{(2)}\}$ with symmetric irreflexive relation $R$ and with
elements $a_1,a_2,a_3$ such that $R=\{[a_1,a_3]$, $[a_2,a_3]\}$.
The substructures $\mathcal{M}_1\subset\mathcal{M}$ and
$\mathcal{M}_2\subset\mathcal{M}$ with the universes $\{a_1\}$ and
$\{a_2\}$ form the substructure $\mathcal{M}(M_1\cup M_2)$ with
the universe $\{a_1,a_2\}$ and it does not satisfy the formula
$\varphi(a_1,a_2)$ in Example \ref{exess3}. At the same time the
structure $\mathcal{M}$ satisfies this formula.}
\end{example}

Since in some cases elementary substructures of given structure
$\mathcal{M}$ form the lattice with respect to the operations
$\mathcal{M}_1\wedge\mathcal{M}_2=\mathcal{M}_1\cap\mathcal{M}_2$
and $\mathcal{M}_1\vee\mathcal{M}_2=\mathcal{M}(M_1\cup M_2)$, the
study of hypergraphs $\mathcal{H}(\mathcal{M})$, for these cases,
is reduced to study of the lattices $L(\mathcal{M})$. As Example
in \cite{Wilkie} shows, the lattices $L(\mathcal{M})$ can be
non-distributive unlike the description in Theorem \ref{thess2},
where correspondent lattices are distributive, and for finite
$H(\mathcal{M}_0)$ even form Boolean algebras.

In the given context hypergraphs/lattices with minimal, i.e. least
structures play an important role. These structures can be
obtained from an arbitrary structure by addition of constants
interpreted by all elements of the structure. Besides, these
minimal structures exist for finite sets $H(\mathcal{M})$.

In \cite{Tan}, the following theorem on dichotomy for minimal
structures is proved.

\begin{theorem}\label{thess3}
Let $\mathcal{M}_0$ be a minimal structure, $\mathcal{M}$ be its
saturated elementary extension and $p\in S_1(\mathcal{M}_0)$ be
unique non-algebraic $1$-type. Then exactly one of the following
conditions holds:

{\rm (I)} the structure $(p(M), {\rm Sem}_p)$ is a pregeometry,
where ${\rm Sem}_p$ is the relation of semi-isolation on the set
of realizations of the type $p$, i.e. the following conditions are
satisfied:

{\rm (S1)} Monotony: if $A\subseteq B$ then $A\subseteq{\rm
Sem}_p(A)\subseteq{\rm Sem}_p(B)$;

{\rm (S2)} Finite character: ${\rm Sem}_p(A)=\bigcup\{{\rm
Sem}_p(A_0)\mid A_0$ is a finite subset of $A\}$;

{\rm (S3)} Transitivity: ${\rm Sem}_p(A)={\rm Sem}_p({\rm
Sem}_p(A))$;

{\rm (S4)} Exchange property {\rm (}Symmetry{\rm )}: if $a\in{\rm
Sem}_p(A\cup\{b\})\setminus{\rm Sem}_p(A)$ then $b\in {\rm
Sem}_p(A\cup\{a\})$;

{\rm (II)} for some finite $A\subset M$ there exists an infinite
set $C_0\subseteq{\rm dcl}(A\cup M_0)$ and a definable quasi-order
$\leq$ on $\mathcal{M}$ such that $C_0$ orders a type over $A$:

{\rm (D1)} for any $c\in C_0$ the set $\{x\in C_0\mid c\leq x\}$
is a cofinite subset of $C_0$;

{\rm (D2)} $C_0$ is an initial segment of $\mathcal{M}$: if $c\in
C_0$ and $m\leq c$, then $m\in C_0$.
\end{theorem}

Basic examples illustrating Theorem \ref{thess3} are represented
by ordered structures $\langle\omega,<\rangle$ and
$\langle\omega+\omega^\ast,<\rangle$. The conclusion of Theorem
\ref{thess2} holds for both structures. Moreover, for
$\mathcal{M}_1\equiv\langle\omega,<\rangle$ and
$\mathcal{M}_2\equiv\langle\omega+\omega^\ast,<\rangle$ the
structures $L(\mathcal{M}_1)$ and $L(\mathcal{M}_2)$ form atomic
Boolean algebras, whose atoms are defined by equivalence classes,
being closures of singletons, not in $\omega+\omega^\ast$, taking
all predecessors and successors.

\medskip
Return to Example \ref{exess2}. It is known that the intersection
of convex sets is convex, whereas the intersection of dense orders
can be not dense. For instance, any interval $[a,b]$ contains
countable dense subsets $X,Y$ such that $X\cap Y=\{a,b\}$. It
means that for the structure $\mathcal{M}'\equiv\langle{\bf
Q},<,c_q\rangle_{q\in{\bf Q}}$ the structure $L(\mathcal{M}')$
forms a lattice, moreover, a Boolean algebra, if and only if each
type in $S_1({\rm Th}(\mathcal{M}'))$ has at most one realization
in $\mathcal{M}'$. If $\mathcal{M}'$, with the lattice
$L(\mathcal{M}')$, realizes $\lambda$ non-principal $1$-types,
then $|L(\mathcal{M}')|=2^\lambda$. Thus, the following
proposition holds.

\begin{proposition}\label{press5}
For the structure $L(\mathcal{M}')$ the following are equivalent:

$(1)$ $L(\mathcal{M}')$ is a lattice;

$(2)$ $L(\mathcal{M}')$ forms an atomic Boolean algebra;

$(3)$ each type in $S_1({\rm Th}(\mathcal{M}'))$ has at most one
realization in $\mathcal{M}'$, and if $\mathcal{M}'$ realizes
$\lambda$ non-principal $1$-types, then
$|L(\mathcal{M}')|=2^\lambda$.
\end{proposition}

Proposition \ref{press5} admits natural modifications for a series
of theories with minimal models, for instance, for models,
obtained by replacement of elements in $\mathcal{M}'$ with finite
antichains of fixed cardinality marked by unary predicates $P_q$
instead of constants $c_q$. Note that admitting replacement of
constants $c_q$ by infinite antichains $P_q$ the structure
$L(\mathcal{M}')$ is not a lattice since $P_q$ can be divided by
some elementary substructures $\mathcal{M}'_1, \mathcal{M}'_2\prec
\mathcal{M}'$ into two disjoint parts, whence $\mathcal{M}'_1\cap
\mathcal{M}'_2\not\prec \mathcal{M}'$.

\medskip
Clearly, as above, in the general case if there are separable
elements in definable sets $A\subseteq M$ of structure
$\mathcal{M}$ then $L(\mathcal{M})$ is not closed under
intersections, i.e., $L(\mathcal{M})$ is not even a lower
semilattice. Thus, the following proposition holds.

\begin{proposition}\label{press6}
If $L(\mathcal{M})$ is a lattice then $\mathcal{M}$ does not have
definable sets $A\subseteq M$ containing elements separable each
other, in particular, $\mathcal{M}$ does not contain
$\mathcal{H}$-free sets $A\subseteq M$.
\end{proposition}

In view of Proposition \ref{press6} it is natural, for given
structure $\mathcal{M}$, along with $L(\mathcal{M})$ to consider
for sets $X\subseteq M$ the following {\em relative} structures
$L_X(\mathcal{M})$. Denote by $H_X(\mathcal{M}$ the family of all
sets in $H(\mathcal{M}$ containing the set $X$. Then
$L_X(\mathcal{M})\rightleftharpoons\langle
H_X(\mathcal{M};\wedge,\vee\rangle$, where for structures
$\mathcal{M}_1,\mathcal{M}_2\prec\mathcal{M}$ containing $X$,
$\mathcal{M}_1\wedge\mathcal{M}_2=\mathcal{M}_1\cap\mathcal{M}_2$
and $\mathcal{M}_1\vee\mathcal{M}_2=\mathcal{M}(M_1\cup M_2)$.

Note that if $X$ is a universe of some elementary substructure of
structure $\mathcal{M}$ then definable sets $A\subseteq M$ already
do not contain elements separable by sets in $L_X(\mathcal{M})$.
Then, in any case, $\mathcal{M}_1\wedge\mathcal{M}_2$ is a
substructure of $\mathcal{M}$ and the elementarity of that
substructure is characterized by Theorem \ref{thTV2}.

\medskip
The following example illustrates that apart from the density
there are other reasons preventing to consider $L(\mathcal{M})$ as
a lattice.

\begin{example}\label{ex_albai_1}\rm \cite{albai} Let $\mathcal{M}=\langle M; <,P^1, U^2, c_i \rangle_{i\in\omega}$ be a linearly
ordered structure such that $\mathcal{M}$ is a disjoint union of
interpretations of unary predicates $P$ and $\neg P$, where $\neg
P(\mathcal{M})<P(\mathcal{M})$. We identify interpretations of $P$
and $\neg P$ with the set $\mathbb{Q}$ of rational numbers with
the natural order.

The symbol $U$ interprets the binary relation defined as follows:
for any $a\in P(\mathcal{M}), b\in \neg P(\mathcal{M})$ $U(a,b)
\Leftrightarrow b<a+\sqrt{2}$.

The constants $c_i$ interpret an infinite strictly increasing
sequence on $P(\mathcal{M})$ as follows: $c_i=i\in \mathbb{Q}$.

Clearly that $Th(\mathcal{M})$ is a weakly o-minimal theory. Let
$$p(x):=\{x>c_i\mid i\in \omega\}\cup \{P(x)\},$$
$$q(y):=\{\forall t(U(c_i, t)\to t<y)\mid i\in\omega\}\cup \{\neg P(y)\}.$$

Obviously, $p, q\in S_1(\emptyset)$ are nonisolated types and
$p\not\perp^w q$. Since there are no $\emptyset$-definable
bijections from $p(\mathcal{M}')$ onto $q(\mathcal{M}')$, where
$\mathcal{M}'$ is a model of $Th(\mathcal{M})$ realizing some of
these types then $Th(\mathcal{M})$ is not quite o-minimal.

As shown in \cite{albai}, $Th(\mathcal{M})$ has exactly 4 pairwise
non-isomorphic countable models: the prime model $\mathcal{M}$,
i.e., with $p(\mathcal{M})=\emptyset$ and
$q(\mathcal{M})=\emptyset$; the model $\mathcal{M}_1$ such that
$p(\mathcal{M}_1)$ has the ordering type $[0,1)\cap\mathbb{Q}$,
$q(\mathcal{M}_1)$ has the ordering type $(0,1)\cap \mathbb{Q}$;
the model $\mathcal{M}_2$ such that $p(M_2)$ has the ordering type
$(0,1)\cap \mathbb{Q}$, $q(M_2)$ has the ordering type
$[0,1)\cap\mathbb{Q}$; and the countable saturated model
$\mathcal{M}_3$.

Therefore $\mathcal{M}_1\cap \mathcal{M}_2\not\prec
\mathcal{M}_3$. By this reason as well as by the possibility of
violation of density in intersections, the structure
$L(\mathcal{M}_3)$ does not form a lower semilattice.
\end{example}

\begin{note}\label{note_ehr_1} \rm
Along with Example \ref{ex_albai_1} if we consider the known
Ehrenfeucht's example with three models: a prime model
$\mathcal{M}_0$, a weakly saturated model $\mathcal{M}_1$, and a
countable saturated model $\mathcal{M}_2$, then the structure
$L(\mathcal{M}_2)$ is not a lattice in view of presence of dence
definable intervals but includes the three-element linearly
ordered lattice consisting of the universes $M_0$, $M_1$, $M_2$.
\end{note}

\section{Lattice structures on sets of isomorphism types of models of a theory}

Following Example \ref{ex_albai_1} and Remark \ref{note_ehr_1} we
consider a question on existence of natural lattices associated
with hypergraphs $(M,H(\mathcal{M}))$ which a distinct to
$L(\mathcal{M})$. Related lattices are lattices represented by
Rudin--Keisler preorders ${\rm RK}(T)$ \cite{CCMCT14} for
isomorphism types of prime models of a theory $T$, over finite
sets, or their lattice fragments.

The description \cite{KulSudRK} of structures ${\rm RK}(T)$ for
Ehrenfeucht quite o-minimal theories $T$ implies that these
structures, for the considered theories, form finite lattices
${\rm LRK}(T)$ consisting of $2^k\cdot 3^s$ elements and, in view
of the main result of the paper \cite{KS}, the number
$I(T,\omega)$ of pairwise non-isomorphic countable models of $T$
equals $3^k\cdot 6^s$, $k,s\in\omega$.

The Hasse diagrams illustrating these lattices ${\rm LRK}(T)$ are
represented in Fig.~\ref{fig1}--\ref{fig9} for the following
values $k$ and $s$:

1) $k=1$, $s=0$;

2) $k=0$, $s=1$;

3) $k=2$, $s=0$;

4) $k=3$, $s=0$;

5) $k=0$, $s=2$;

6) $k=0$, $s=3$;

7) $k=1$, $s=1$;

3) $k=2$, $s=1$;

5) $k=1$, $s=2$.

\begin{figure}[t]
\begin{center}
\unitlength 4mm
\begin{picture}(5,6.5)(-5.5,1.0)
{\footnotesize\put(2,2.5){\line(0,5){5}}
\put(2,2.5){\makebox(0,0)[cc]{$\bullet$}}
\put(2,7.5){\makebox(0,0)[cc]{$\bullet$}} }
\end{picture}
\hfill \unitlength 5mm
\begin{picture}(22,7)(4,1.7)
{\footnotesize \put(19.5,2.5){\line(0,5){5}}
\put(19.5,2.5){\makebox(0,0)[cc]{$\bullet$}}
\put(19.5,5){\makebox(0,0)[cc]{$\bullet$}}
\put(19.5,7.5){\makebox(0,0)[cc]{$\bullet$}}
 }
\end{picture}
\end{center}
\parbox[t]{0.4\textwidth}{\caption{}\label{fig1}}
\hfill
\parbox[t]{0.4\textwidth}{\caption{}\label{fig2}}

\end{figure}

\begin{figure}[t]
\begin{center}
\unitlength 14mm
\begin{picture}(5,1)(-1.14,-0.1)
{\footnotesize \put(0,1){\makebox(0,0)[cc]{$\bullet$}}

\put(1,0){\makebox(0,0)[cc]{$\bullet$}}
\put(2,1){\makebox(0,0)[cc]{$\bullet$}}

\put(1,2){\makebox(0,0)[cc]{$\bullet$}} \put(1,0){\line(1,1){1}}
\put(1,0){\line(-1,1){1}} \put(1,2){\line(-1,-1){1}}
\put(1,2){\line(1,-1){1}}
 }
\end{picture}
\hfill \unitlength 14mm
\begin{picture}(5,3.5)(1.35,0.2)
{\footnotesize\put(3,1){\makebox(0,0)[cc]{$\bullet$}}
\put(3,2){\makebox(0,0)[cc]{$\bullet$}}
\put(4,0){\makebox(0,0)[cc]{$\bullet$}}
\put(4,1){\makebox(0,0)[cc]{$\bullet$}}
\put(4,2){\makebox(0,0)[cc]{$\bullet$}}
\put(4,3){\makebox(0,0)[cc]{$\bullet$}}
\put(5,1){\makebox(0,0)[cc]{$\bullet$}}
\put(5,2){\makebox(0,0)[cc]{$\bullet$}} \put(3,1){\line(0,1){1}}
\put(5,1){\line(0,1){1}} \put(4,0){\line(0,1){1}}
\put(4,2){\line(0,1){1}} \put(4,0){\line(1,1){1}}
\put(4,1){\line(1,1){1}} \put(4,0){\line(-1,1){1}}
\put(4,1){\line(-1,1){1}} \put(4,3){\line(-1,-1){1}}
\put(4,2){\line(-1,-1){1}} \put(4,3){\line(1,-1){1}}
\put(4,2){\line(1,-1){1}} \put(4.2,0){\makebox(0,0)[cl]{0}} }
\end{picture}
\end{center}
\parbox[t]{0.4\textwidth}{\caption{}\label{fig3}}
\hfill
\parbox[t]{0.4\textwidth}{\caption{}\label{fig4}}

\end{figure}

\begin{figure}
\begin{center}
\unitlength 18mm
\begin{picture}(3,1)(-0.69,-0.1)
{\footnotesize \put(0,1){\makebox(0,0)[cc]{$\bullet$}}

\put(1,0){\makebox(0,0)[cc]{$\bullet$}}
\put(2,1){\makebox(0,0)[cc]{$\bullet$}}

\put(1,2){\makebox(0,0)[cc]{$\bullet$}}
\put(1,0){\line(1,1){1}}\put(1.5,0.5){\line(-1,1){1}}
\put(0.5,0.5){\line(1,1){1}} \put(1,0){\line(-1,1){1}}
\put(1,2){\line(-1,-1){1}} \put(1,2){\line(1,-1){1}}
\put(1.5,0.5){\makebox(0,0)[cc]{$\bullet$}}
\put(0.5,0.5){\makebox(0,0)[cc]{$\bullet$}}
\put(1,1){\makebox(0,0)[cc]{$\bullet$}}
\put(0.5,1.5){\makebox(0,0)[cc]{$\bullet$}}
\put(1.5,1.5){\makebox(0,0)[cc]{$\bullet$}}
 }
\end{picture}
\hfill \unitlength 20mm
\begin{picture}(4,3.5)(1.6,0.2)
{\footnotesize\put(3,1){\makebox(0,0)[cc]{$\bullet$}}
\put(3,2){\makebox(0,0)[cc]{$\bullet$}}
\put(4,0){\makebox(0,0)[cc]{$\bullet$}}
\put(4,1){\makebox(0,0)[cc]{$\bullet$}}
\put(4.006,2){\makebox(0,0)[cc]{$\bullet$}}
\put(4,3){\makebox(0,0)[cc]{$\bullet$}}
\put(5,1){\makebox(0,0)[cc]{$\bullet$}}
\put(5,2){\makebox(0,0)[cc]{$\bullet$}} \put(3,1){\line(0,1){1}}
\put(5,1){\line(0,1){1}} \put(4,0){\line(0,1){1}}
\put(4,2){\line(0,1){1}} \put(4,0){\line(1,1){1}}
\put(4,1){\line(1,1){1}} \put(4,0){\line(-1,1){1}}
\put(4,1){\line(-1,1){1}} \put(4,3){\line(-1,-1){1}}
\put(4,2){\line(-1,-1){1}} \put(4,3){\line(1,-1){1}}
\put(4,2){\line(1,-1){1}}
\put(4,2.5){\makebox(0,0)[cc]{$\bullet$}}
\put(3.35,0.65){\makebox(0,0)[cc]{$\bullet$}}
\put(3.35,1.15){\makebox(0,0)[cc]{$\bullet$}}
\put(3.35,1.65){\makebox(0,0)[cc]{$\bullet$}}
\put(3.95,1.75){\makebox(0,0)[cc]{$\bullet$}}
\put(3.95,1.25){\makebox(0,0)[cc]{$\bullet$}}
\put(3.945,2.25){\makebox(0,0)[cc]{$\bullet$}}
\put(4.35,1.65){\makebox(0,0)[cc]{$\bullet$}}
\put(4.35,2.15){\makebox(0,0)[cc]{$\bullet$}}
\put(4.35,2.65){\makebox(0,0)[cc]{$\bullet$}}
\put(3.6,1.6){\makebox(0,0)[cc]{$\bullet$}}
\put(3.6,2.1){\makebox(0,0)[cc]{$\bullet$}}
\put(3.6,2.6){\makebox(0,0)[cc]{$\bullet$}}
\put(4.6,0.6){\makebox(0,0)[cc]{$\bullet$}}
\put(4.6,1.1){\makebox(0,0)[cc]{$\bullet$}}
\put(4.6,1.6){\makebox(0,0)[cc]{$\bullet$}}
\put(4,0.5){\makebox(0,0)[cc]{$\bullet$}}
\put(5,1.5){\makebox(0,0)[cc]{$\bullet$}}
\put(3,1.5){\makebox(0,0)[cc]{$\bullet$}}

\put(3.35,0.65){\line(0,1){1}} \put(4.6,0.6){\line(0,1){1}}
\put(3.35,0.65){\line(1,1){1}} \put(4.6,0.6){\line(-1,1){1}}
\put(4,0.5){\line(1,1){1}} \put(4,0.5){\line(-1,1){1}}
\put(3.35,1.65){\line(1,1){1}} \put(4.6,1.6){\line(-1,1){1}}
\put(3,1.5){\line(1,1){1}} \put(5,1.5){\line(-1,1){1}}
\put(3.6,1.6){\line(0,1){1}} \put(4.35,1.65){\line(0,1){1}}
\put(3.35,1.15){\line(1,1){1}} \put(4.6,1.1){\line(-1,1){1}}
\put(3.95,1.25){\line(0,1){1}}

}
\end{picture}
\end{center}
\parbox[t]{0.4\textwidth}{\caption{}\label{fig5}}
\hfill
\parbox[t]{0.4\textwidth}{\caption{}\label{fig6}}

\end{figure}

\begin{figure}[t]
\begin{center}
\unitlength 24mm
\begin{picture}(3,1)(-0.05,-0.1)
{\footnotesize \put(1,0){\makebox(0,0)[cc]{$\bullet$}}
\put(2,1){\makebox(0,0)[cc]{$\bullet$}}
\put(1,0){\line(1,1){1}}\put(1.5,0.5){\line(-1,1){0.5}}
\put(0.5,0.5){\line(1,1){1}} \put(1,0){\line(-1,1){0.5}}
\put(1.5,1.5){\line(1,-1){0.5}}
\put(1.5,0.5){\makebox(0,0)[cc]{$\bullet$}}
\put(0.5,0.5){\makebox(0,0)[cc]{$\bullet$}}
\put(1,1){\makebox(0,0)[cc]{$\bullet$}}

\put(1.5,1.5){\makebox(0,0)[cc]{$\bullet$}}

 }
\end{picture}
\hfill \unitlength 24mm
\begin{picture}(3,2.5)(2.65,0.1)
{\footnotesize

\put(4,0){\makebox(0,0)[cc]{$\bullet$}}
\put(4,1){\makebox(0,0)[cc]{$\bullet$}}

\put(5,1){\makebox(0,0)[cc]{$\bullet$}}

\put(5,1){\line(0,1){0.5}} \put(4,0){\line(0,1){0.5}}
 \put(4,0){\line(1,1){1}}
\put(4,0){\line(-1,1){0.5}}

 \put(4.5,1.5){\line(1,-1){0.5}}

\put(3.5,0.5){\makebox(0,0)[cc]{$\bullet$}}
\put(3.5,1.0){\makebox(0,0)[cc]{$\bullet$}}

\put(4.0,1.5){\makebox(0,0)[cc]{$\bullet$}}
\put(4.5,1.5){\makebox(0,0)[cc]{$\bullet$}}
\put(4.5,2.0){\makebox(0,0)[cc]{$\bullet$}}

\put(4.5,0.5){\makebox(0,0)[cc]{$\bullet$}}
\put(4.5,1){\makebox(0,0)[cc]{$\bullet$}}

\put(4,0.5){\makebox(0,0)[cc]{$\bullet$}}
\put(5,1.5){\makebox(0,0)[cc]{$\bullet$}}

\put(3.5,0.5){\line(0,1){0.5}} \put(4.5,0.5){\line(0,1){0.5}}
\put(3.5,0.5){\line(1,1){1}}
\put(4.5,0.5){\line(-1,1){0.5}}\put(4.5,1){\line(-1,1){0.5}}
\put(4,0.5){\line(1,1){1}} \put(4,0.5){\line(-1,1){0.5}}

\put(5,1.5){\line(-1,1){0.5}} \put(4.5,1.5){\line(0,1){0.5}}
\put(3.5,1.0){\line(1,1){1}} \put(4,1){\line(0,1){0.5}}

}
\end{picture}
\end{center}
\parbox[t]{0.4\textwidth}{\caption{}\label{fig7}}
\hfill
\parbox[t]{0.4\textwidth}{\caption{}\label{fig8}}

\end{figure}

\begin{figure}[t]
\begin{center}
\unitlength 21mm
\begin{picture}(4,2.95)(2.05,0.0)
{\footnotesize\put(3,1){\makebox(0,0)[cc]{$\bullet$}}

\put(4,0){\makebox(0,0)[cc]{$\bullet$}}
\put(4,1){\makebox(0,0)[cc]{$\bullet$}}
\put(4.0,2){\makebox(0,0)[cc]{$\bullet$}}

\put(5,1){\makebox(0,0)[cc]{$\bullet$}}

\put(3,1){\line(0,1){0.5}} \put(5,1){\line(0,1){0.5}}
\put(4,0){\line(0,1){0.5}} \put(4,2){\line(0,1){0.5}}
\put(4,0){\line(1,1){1}} \put(4,0){\line(-1,1){1}}
 \put(4,2){\line(-1,-1){1}}
 \put(4,2){\line(1,-1){1}}
\put(4,2.5){\makebox(0,0)[cc]{$\bullet$}}
\put(3.5,0.5){\makebox(0,0)[cc]{$\bullet$}}
\put(3.5,1.0){\makebox(0,0)[cc]{$\bullet$}}
\put(3.5,1.5){\makebox(0,0)[cc]{$\bullet$}}
\put(3.5,2.0){\makebox(0,0)[cc]{$\bullet$}}
\put(4.0,1.5){\makebox(0,0)[cc]{$\bullet$}}
\put(4.5,1.5){\makebox(0,0)[cc]{$\bullet$}}
\put(4.5,2.0){\makebox(0,0)[cc]{$\bullet$}}

\put(4.5,0.5){\makebox(0,0)[cc]{$\bullet$}}
\put(4.5,1){\makebox(0,0)[cc]{$\bullet$}}

\put(4,0.5){\makebox(0,0)[cc]{$\bullet$}}
\put(5,1.5){\makebox(0,0)[cc]{$\bullet$}}
\put(3,1.5){\makebox(0,0)[cc]{$\bullet$}}

\put(3.5,0.5){\line(0,1){0.5}} \put(4.5,0.5){\line(0,1){0.5}}
\put(3.5,0.5){\line(1,1){1}}
\put(4.5,0.5){\line(-1,1){1}}\put(4.5,1){\line(-1,1){1}}
\put(4,0.5){\line(1,1){1}} \put(4,0.5){\line(-1,1){1}}
 \put(3,1.5){\line(1,1){1}}
\put(5,1.5){\line(-1,1){1}} \put(3.5,1.5){\line(0,1){0.5}}
\put(4.5,1.5){\line(0,1){0.5}} \put(3.5,1.0){\line(1,1){1}}
\put(4,1){\line(0,1){0.5}}

}
\end{picture}
\end{center}
{ \caption{} \label{fig9}}

\end{figure}

\begin{theorem}\label{th_qom_l} Let $T$ be an Ehrenfeucht quite
o-minimal theory, $I(T,\omega)=3^k\cdot 6^s$, $k,s\in\omega$.
Then:

$(1)$ ${\rm LRK}(T)$ is a lattice;

$(2)$ ${\rm LRK}(T)$ is a Boolean algebra $\Leftrightarrow$ $k\geq
1$ and $s=0$; in such a case the Boolean lattice ${\rm LRK}(T)$
has a cardinality $2^k$;

$(3)$ ${\rm LRK}(T)$ is linearly ordered $\Leftrightarrow$
$k+s\leq 1$.
\end{theorem}

Proof of Theorem \ref{th_qom_l}. Let $\Gamma=\Gamma_1\cup\Gamma_2$
be a maximal independent set of nonisolated types in $S_1(T)$,
where realizations of each type in $\Gamma_1$ generate three
models, with prime one, and realizations of each type in
$\Gamma_2$ generate six models, with prime one, $|\Gamma_1|=k$,
$|\Gamma_2|=s$.

$(1)$ We argue to show that ${\rm LRK}(T)$ is a lattice. Indeed,
for isomorphism types $\widetilde{\mathcal{M}_1}$ and
$\widetilde{\mathcal{M}_2}$ of prime model $\mathcal{M}_1$ and
$\mathcal{M}_2$ over some finite sets $A$ and $B$, respectively,
we define sets $X,Y\subseteq\Gamma\times\{0,1\}$ defining these
isomorphism types such that $X=\{(p,0)\mid \mathcal{M}_1\models
p(a)\mbox{ for some }a\in A,\mbox{ and }|p(\mathcal{M}_1)|=1\mbox{
or }p\in\Gamma_1\}\cup \{(p,1)\mid \mathcal{M}_1\models p(a)\mbox{
for some }a\in A,|p(\mathcal{M}_1)|\geq\omega,p\in\Gamma_2\}$ and
$Y=\{(q,0)\mid \mathcal{M}_2\models q(b)\mbox{ for some }b\in
B,\mbox{ and }|q(\mathcal{M}_2)|=1\mbox{ or }q\in\Gamma_1\}\cup
\{(q,1)\mid \mathcal{M}_2\models q(b)\mbox{ for some }b\in
B,|q(\mathcal{M}_2)|\geq\omega,q\in\Gamma_2\}$. Then the
isomorphism type for
$\widetilde{\mathcal{M}_1}\wedge\widetilde{\mathcal{M}_2}$
corresponds to the set $U\subseteq\Gamma\times\{0,1\}$ consisting
of all common pairs of $X$ and $Y$, as well as all possible pairs
$(p,0)$, if $(p,0)\in X$ and $(p,1)\in Y$, or $(p,1)\in X$ and
$(p,0)\in Y$. And the isomorphism type for
$\widetilde{\mathcal{M}_1}\vee\widetilde{\mathcal{M}_2}$
corresponds to the set $V\subseteq\Gamma\times\{0,1\}$ consisting
of the following pairs:

i) all common pairs of $X$ and $Y$,

ii) all pairs $(p,i)\in X$ such that
$Y\cap\{(p,0),(p,1)\}\emptyset$,

iii) all pairs $(p,i)\in Y$ such that
$X\cap\{(p,0),(p,1)\}\emptyset$,

iv) all pairs $(p,1)$ such that $(p,0)\in X$ and $(p,1)\in Y$, or
$(p,1)\in X$ and $(p,0)\in Y$.

The defined correspondence witnesses that ${\rm LRK}(T)$ is a
lattice.

$(2)$ If $s\ne 0$ then ${\rm LRK}(T)$ is not a Boolean algebra by
Stone Theorem, since the cardinality of each finite Boolean
algebra equals $2^n$ for some $n\in\omega$ whereas $|{\rm
LRK}(T)|=2^k\cdot 3^s$. If $s=0$ then ${\rm LRK}(T)$ is a Boolean
algebra of a cardinality $2^k$ such that for isomorphism types
$\widetilde{\mathcal{M}_1}$ and $\widetilde{\mathcal{M}_2}$ of
prime models $\mathcal{M}_1$ and $\mathcal{M}_2$ over some finite
sets $A$ and $B$, respectively, and for sets $X,Y\subseteq\Gamma$
such that $X=\{p(x)\in\Gamma\mid \mathcal{M}_1\models p(a)\mbox{
for some }a\in A \}$ and $Y=\{q(x)\in\Gamma\mid
\mathcal{M}_2\models q(b)\mbox{ for some }b\in B \}$, the
isomorphism type
$\widetilde{\mathcal{M}_1}\wedge\widetilde{\mathcal{M}_2}$
corresponds to the set $X\cap Y$, and the isomorphism type
$\widetilde{\mathcal{M}_1}\vee\widetilde{\mathcal{M}_2}$
corresponds to the set $X\cup Y$.

(3) If $k+s\leq 1$ then ${\rm LRK}(T)$ is linearly ordered as
shown in Fig. \ref{fig1} and \ref{fig2}. If $k+s
>1$ then $|\Gamma|>1$ and for distinct types $p,q\in\Gamma$ the isomorphism types
of models $\mathcal{M}_p$ and $\mathcal{M}_q$ are incomparable in
${\rm LRK}(T)$. \endproof

\medskip
The description for distributions of disjoint unions of
Ehrenfeucht theories and the arguments for the proof of Theorem
\ref{th_qom_l} allow to formulate the following theorem modifying
Theorem  \ref{th_qom_l}.

\begin{theorem}\label{th_2} Let $T$ be a disjoint union of theories $T_1$ and $T_2$
in disjoint languages and having finite numbers $I(T_1,\omega)$
and $I(T_2,\omega)$ of countable models. Then:

$(1)$ ${\rm LRK}(T)$ is a {\rm (}Boolean{\rm )} lattice
$\Leftrightarrow$ ${\rm LRK}(T_1)$ and ${\rm LRK}(T_2)$ are {\rm
(}Boolean{\rm )} lattices;

$(2)$ ${\rm LRK}(T)$ is linearly ordered $\Leftrightarrow$  ${\rm
LRK}(T_1)$ and ${\rm LRK}(T_2)$ are linearly ordered, \ and \
${\rm min}\{I(T_1,\omega),I(T_2,\omega)\}=1$.
\end{theorem}

Proof. (1) If ${\rm LRK}(T)$ is a {\rm (}Boolean{\rm )} lattice,
then ${\rm LRK}(T_1)$ and ${\rm LRK}(T_2)$ are {\rm (}Boolean{\rm
)} lattices, since ${\rm LRK}(T_1)$ and ${\rm LRK}(T_2)$ are
isomorphic to sublattices $L_1$ and $L_2$ of the lattice ${\rm
LRK}(T)$, and elements/complements of elements in ${\rm LRK}(T)$
are represented as pairs of elements/complements of elements in
$L_1$ and $L_2$. If ${\rm LRK}(T_1)$ and ${\rm LRK}(T_2)$ are {\rm
(}Boolean{\rm )} lattices, then ${\rm LRK}(T)$ is a {\rm
(}Boolean{\rm )} lattice again in view of aforesaid
representation.

(2) If ${\rm LRK}(T)$ is linearly ordered then ${\rm LRK}(T_1)$
and ${\rm LRK}(T_2)$ are linearly ordered, being isomorphic to
substructures of ${\rm LRK}(T)$. Here $T_1$ or $T_2$ should be
$\omega$-categorical, since otherwise prime models over pairs
$(p_1,q_1)$ and $(p_2,q_2)$ occur ${\rm LRK}(T)$-incomparable,
where $p_1,p_2\in S_1(T_1)$, $q_1,q_2\in S_1(T_2)$, $p_1,q_2$ are
isolated, $p_2,q_1$ are nonisolated.

If structures ${\rm LRK}(T_1)$ and ${\rm LRK}(T_2)$ linearly
ordered, and ${\rm min}\{I(T_1,\omega),I(T_2,\omega)\}=1$, then
${\rm LRK}(T)$ is linearly ordered, since ${\rm LRK}(T)\simeq {\rm
LRK}(T_1)$ for $I(T_2,\omega)=1$, and ${\rm LRK}(T)\simeq {\rm
LRK}(T_2)$ for $I(T_1,\omega)=1$.
\endproof

\medskip
In Fig.~\ref{fig10} and \ref{fig11} we illustrate Theorem
\ref{th_2} by structures ${\rm LRK}(T)$ in \cite{SudRK}, for
disjoint unions of theories, which are not lattices.

\begin{figure}
\vspace{5mm}
\begin{center}
\unitlength 4mm
\begin{picture}(7,5)(25.5,8.5)
{\footnotesize

\put(14.5,2.5){\line(-2,5){2}} \put(14.5,2.5){\line(2,5){2}}
\put(14.5,2.5){\makebox(0,0)[cc]{$\bullet$}}
\put(12.5,7.5){\makebox(0,0)[cc]{$\bullet$}}
\put(16.5,7.5){\makebox(0,0)[cc]{$\bullet$}}
\put(14.5,2.5){\circle{0.6}} \put(14.5,7.5){\oval(5,1)}
%\put(15.5,2.5){\makebox(0,0)[cr]{$0$}}
%\put(11.7,7.5){\makebox(0,0)[cr]{$k$}}

\put(20.5,6.5){\line(-2,5){2}} \put(20.5,6.5){\line(2,5){2}}
\put(20.5,6.5){\makebox(0,0)[cc]{$\bullet$}}
\put(18.5,11.5){\makebox(0,0)[cc]{$\bullet$}}
\put(22.5,11.5){\makebox(0,0)[cc]{$\bullet$}}
\put(20.5,6.5){\circle{0.6}} \put(20.5,11.5){\oval(5,1)}
%\put(21.7,6.5){\makebox(0,0)[cr]{$m$}}
%\put(22.9,12.5){\makebox(0,0)[cr]{$k+2m+km$}}

\put(14.5,2.5){\line(3,2){6}} \put(12.5,7.5){\line(3,2){6}}
\put(16.5,7.5){\line(3,2){6}}

\put(12.5,7.5){\line(5,2){10}} \put(16.5,7.5){\line(2,4){2}}

}
\end{picture}

\begin{picture}(7,5)(7.5,3.5)
{\footnotesize

\put(14.5,2.5){\line(-2,5){2}} \put(14.5,2.5){\line(2,5){2}}
\put(14.5,2.5){\makebox(0,0)[cc]{$\bullet$}}
\put(12.5,7.5){\makebox(0,0)[cc]{$\bullet$}}
\put(16.5,7.5){\makebox(0,0)[cc]{$\bullet$}}
\put(14.5,2.5){\circle{0.6}} \put(14.5,7.5){\oval(5,1)}
%\put(15.5,2.5){\makebox(0,0)[cr]{$0$}}
%\put(11.7,7.5){\makebox(0,0)[cr]{$k$}}

\put(20.5,6.5){\line(-2,5){2}} \put(20.5,6.5){\line(2,5){2}}
\put(20.5,6.5){\makebox(0,0)[cc]{$\bullet$}}
\put(18.5,11.5){\makebox(0,0)[cc]{$\bullet$}}
\put(22.5,11.5){\makebox(0,0)[cc]{$\bullet$}}
%\put(20.5,6.5){\circle{0.6}}
\put(20.5,11.5){\oval(5,1)}
%\put(21.7,6.5){\makebox(0,0)[cr]{$m$}}
%\put(22.9,12.5){\makebox(0,0)[cr]{$k+2m+km$}}

\put(26.5,6.5){\line(-2,5){2}} \put(26.5,6.5){\line(2,5){2}}
\put(26.5,6.5){\makebox(0,0)[cc]{$\bullet$}}
\put(24.5,11.5){\makebox(0,0)[cc]{$\bullet$}}
\put(28.5,11.5){\makebox(0,0)[cc]{$\bullet$}}
%\put(26.5,6.5){\circle{0.6}}
\put(26.5,11.5){\oval(5,1)}
%\put(27.85,6.5){\makebox(0,0)[cr]{$m$}}
%\put(26.1,12.8){\makebox(0,0)[cr]{$2k+2m+km$}}

\put(26.5,6.5){\line(-2,5){2}} \put(26.5,6.5){\line(2,5){2}}

\put(14.5,2.5){\line(3,2){6}} \put(12.5,7.5){\line(3,2){6}}
\put(16.5,7.5){\line(3,2){6}}

\put(12.5,7.5){\line(4,1){16}}

\put(14.5,2.5){\line(6,2){12}}

\put(12.5,7.5){\line(6,2){12}} \put(16.5,7.5){\line(6,2){12}}

\put(12.5,7.5){\line(5,2){10}} \put(16.5,7.5){\line(4,2){8}}

\put(16.5,7.5){\line(2,4){2}}

\put(23.45,6.5){\oval(7,1)}

\put(23.5,11.5){\oval(12,1.5)}

}
\end{picture}

\end{center}
\parbox[t]{0.4\textwidth}{\caption{}\label{fig10}}
\hfill
\parbox[t]{0.4\textwidth}{\caption{}\label{fig11}}
\end{figure}

\end{document}